\documentclass[letterpaper, 10 pt, conference]{ieeeconf}  

\IEEEoverridecommandlockouts                              
\usepackage{color,xcolor}
\usepackage{ulem}
\usepackage{amsmath}
\newcommand{\sgn}{\mathrm{sgn}}
\usepackage{amsmath,amssymb}
\usepackage{graphicx}
\usepackage{version}
\usepackage{enumerate}
\usepackage{subcaption}
\usepackage{cite}

\newtheorem{Lemma}{Lemma}
\newtheorem{Remark}{Remark}
\newtheorem{thm}{Theorem}
\newtheorem{Corollary}{Corollary}
\newcommand{\MG}[1]{\textcolor{black}{#1}}
\newcommand{\JPB}[1]{\textcolor{black}{#1}}

\makeatother

\newenvironment{figurehere*}
  {\def\@captype{figure*}}
  {}

\title{\LARGE \bf  A Power Tower Control: A New Sliding Mode Control} 

\author{Malek GHANES and Jean-Pierre BARBOT 
\thanks{M. GHANES and J.-P. BARBOT are with Nantes Universite-\'Ecole Centrale de Nantes-LS2N, UMR 6004 CNRS, Nantes, France {\tt\small  malek.ghanes@ec-nantes.fr}}
\thanks{J.-P. BARBOT is also with QUARTZ EA 7393, ENSEA, Cergy-Pontoise, France {\tt\small barbot@ensea.fr}}
}
\begin{document}
\maketitle
\thispagestyle{empty}
\pagestyle{empty}
\begin{abstract}
    \MG{A control based power tower function at order 2 is proposed in this paper. This leads to a new sliding mode control, which allows employing backstepping technique that combines both guaranteed and finite time convergence. The proposed control is applied to a double integrator subject to perturbation $d$. Both guaranteed and finite convergence are ensured by the controller when $d$ is considered constant and bounded, without knowing its upper bound. For the case, when $d$ is variable and bounded with its upper bound known, only a finite time convergence is obtained. Simulation results are given to show the well founded of the proposed novel control.}
\end{abstract}
\begin{keywords} Power tower, Backstepping, Guaranteed-finite time convergence, sliding mode.
\end{keywords}
\section{Introduction}
\MG{There are a large number of iterative techniques to built control or differentiator. These include: High order sliding mode  \cite{levant93, Fridman2002higher,UTKIN2020}, homogeneous \cite{PEP14,Rosier1992homogeneous,Polyakov_book,Andrieu2008homogeneous} backstepping \cite{FREEMAN1993431,KKK_92,Backstepping_Arezki}, singular perturbation \cite{kokotovic1999singular,Vasileva78,Tikhonov52}, high gain \cite{Bornard91,khalil2014high,Khalil2002nonlinear}, ... In this article we study the possibility of constructing a finite-time control using the power tower function \cite{Knoebel81,lynch2017fractal}  truncated to second order. Our original motivation was, similar to the case of the variable exponent “Homogeneous” differentiator, to propose a continuous variable exponent “Homogeneous” control. In the case of the differentiator the  variation law of the homogeneity exponent is a function of the measurement noise \cite{ghanes2020new} or of time to ensure guaranteed and finite convergence \cite{Ghanes2022arbitrary}.  For the control, in order to ensure both guaranteed and finite time convergence with a continuous law exponential variation, theoretical obstructions prevented us from finding the control (see \cite{Wang2023} for a discontinuous law) without a specific continuous power function. }
\MG{Indeed, for our best knowledge, a control combining guaranteed and finite time without singular problems has not yet been considered in the literature. Moreover, this kind of control has not been formally associated with the backstepping approach.}
\MG{ To solve the problem of ensuring at the same time, guaranteed and finite time convergence, the addition of at least two different controllers are proposed in the literature (see \cite{Cruz2017homogeneous}, \cite{Polyakov_book,polyakov2015finite,Overview_Finite_Fixed},...). In this paper we propose to fix this control  problem  by using only one controller. For this purpose, a new control based on the power tower function is introduced. This function, beyond its specific properties at the limits or on the fractal  topology \cite{lynch2017fractal} obtained the property of its derivative, allowed us to render possible the use of backstepping techniques for convergence with both guaranteed\footnote{The convergence is ensured to reach a vicinity of the equilibrium point in a guaranteed time whatever the initial conditions are.} and finite time. By doing so, two parameters, power exponent and linear gain, are needed for the control tuning in absence of perturbation. 
In the presence of the latter,  when it is variable, bounded with its upper bound known, one exponet of the power tower function is set to zero. When the perturbation is constant and bounded, an integral action is necessary.\\
The remaining of the paper is organized as follows. Section \ref{main result} presents the main result, which render possible the use of backstopping technique that ensures a guaranteed-finite convergence by using only one controller. In Section  \ref{simu},  the performances of the proposed control based power tower function, in absence and presence of perturbation, are put forward. A conclusion is given in section \ref{conc} with some future works.} 

\section{Main result} \label{main result}
\MG{Our control exploits the \textbf{power tower function} of order 2 allowing us to propose 
a new \MG{guaranteed}-finite time backstepping based control. Before presenting this new idea,  we need to introduce the following lemma:\\}
\begin{Lemma} \label{lemma}
The time derivative of the function $\lceil a \rfloor^{|a|^\alpha}$ with $a$ a function of time at least $\mathbb{C}^1$ and $\alpha$ a strictly positive constant is \MG{for $a\neq0$}:
\begin{eqnarray}
\frac{d \lceil a \rfloor^{|a|^\alpha}}{d t}&=& |a|^{|a|^\alpha} |a|^{\alpha-1} (1+\alpha\, ln(|a|)) \dot{a} \label{eq:lemme} \\
&=&|a|^{|a|^\alpha+\alpha-1} (1+\alpha\, ln(|a|)) \dot{a} \nonumber
\end{eqnarray}  
with $\dot{a}=\frac{d\, a}{d \, t}$. \vspace*{5mm}
\end{Lemma}
\textbf{Proof}: As 
\[
\lceil a \rfloor^{|a|^\alpha}=| a |^{|a|^\alpha} \sgn(a),
\]
we first time derive 
\[
| a |^{|a|^\alpha}.
\]

\noindent In order to use the derivative of the composition of functions, we set 

\[
| a |^{|a|^\alpha}= e^{ln(|a|) |a|^\alpha}
\]

\noindent and we obtain

\[
\frac{d\,| a |^{|a|^\alpha}}{d t}=e^{ln(|a|) |a|^\alpha} \frac{d \,ln(|a|) |a|^\alpha}{d\, a} \dot{a}
\]

\noindent or again

\[
\frac{d\,| a |^{|a|^\alpha}}{d t}=e^{ln(|a|) |a|^\alpha} (|a|^{\alpha-1}+ln(|a|) \, \alpha |a|^{\alpha-1}) \sgn(a)\dot{a}
\]

\noindent which gives, multiplying both parts by $\sgn(a)$

\[
\frac{d\,\lceil a \rfloor^{|a|^\alpha}}{d t}=e^{ln(|a|) |a|^\alpha} (|a|^{\alpha-1}+ln(|a|) \, \alpha |a|^{\alpha-1})  \dot{a}.
\]

\noindent This ends the proof.\hfill $\Box$\\

\begin{Remark} \label{remark1}
If $\alpha > 1$ then we have 
\begin{equation}
  Lim_{a\rightarrow 0} |a|^{|a|^\alpha+\alpha-1} (1+\alpha\, ln(|a|)) = 0
\end{equation}
which will be a guarantee of a bounded control law and therefore feasible in the vicinity of $a=0$. \MG{This also ensures that for $\alpha>1$ and $\dot{a}$ bounded the equation \eqref{eq:lemme} is also defined at $a=0$ and is equal to zero.}
\end{Remark}
\begin{Remark}\label{remark2}
    \JPB{An other property for $\alpha>0$ of the proposed truncated power tower function is
\begin{equation}
  Lim_{a\rightarrow 0} |a|^{|a|^\alpha} \sgn(a) = Lim_{a\rightarrow 0} \sgn(a). \label{eq:limPT_a=0}
\end{equation}
Consequently the closed loop behavior obtained with such function refers to a sliding mode behavior.}
\end{Remark}

Now, let us consider the following system:
\begin{eqnarray}
    \dot{x}_1&=&x_2 \label{eq:system_ref}  \\
    \dot{x}_2&=&u+d \nonumber
\end{eqnarray}
where $u$ is the control input and $d$ the disturbance.\\

\begin{Remark} \label{Remark3}
    To design a backstepping type control law in \MG{guaranteed}-finite time we will refer to the lemma \ref{lemma} and not to the derivative of $|a|^\theta$ which has an unbounded limit for $1> \theta > 0$ and $a \rightarrow 0^+$.
\end{Remark} 

\begin{thm} \label{Theorem}
    If $d=0$, the following \MG{power tower} control law:
    \begin{eqnarray} \label{pt2}
       u&=&-K_2\lceil z_2 \rfloor^{|z_2|^\gamma}- \left(\frac{}{} x_1  K_1|x_1|^{|x_1|^\beta+\beta-1} \label{u} \right.\\
       &+&\left. (1+\beta\, ln(|x_1|)) (-K_1\lceil x_1 \rfloor^{|x_1|^\beta}+z_2) \frac{}{}\right) \nonumber
    \end{eqnarray}
    where $\beta>1$, $\gamma >0$, $K_1>0$, $K_2>0$ and 
\begin{equation} \label{z_2}
z_2=x_2- x_2^*
\end{equation} 
with $x_2^*= -K_1\lceil x_1 \rfloor^{|x_1|^\beta}$, ensures a \MG{guaranteed}-finite time convergence of \eqref{eq:system_ref} to $x_1=x_2=0$. 
\end{thm}
\vspace*{5mm}
\textbf{Proof}: Based on \MG{the} well known backstepping method \cite{KKK95}, \MG{the proof is given in two steps.} \vspace*{3mm}\\
\MG{\textbf{First step.} 
The first step consists on stabilizing $x_1$ by a fictive control $x_2^*$. For that, we define this control as a power tower control, that is 
\begin{equation} \label{ptx1}
  x_2^*=-K_1\lceil x_1 \rfloor^{|x_1|^\beta}  
\end{equation}
with $\beta >1$. Setting 
\begin{equation} \label{v1}
  V_1=\frac{x_1^2}{2}.  
\end{equation}
The time-derivative of $\eqref{v1}$ reads 
\begin{equation}
    \dot{V}_1=- K_1|x_1|^{1+|x_1|^\beta} + x_1 z_2 \label{dot_V_2}
\end{equation}
\noindent If $x_2^*=x_2$ (i.e., $z_2=0$), then the \MG{guaranteed}-finite time convergence of $x_1$ to zero is achieved when $d=0$. However, at this step, \eqref{z_2} is not converged to zero, that is why the following second step is important to design the real control.}\vspace*{3mm}\\
\MG{\textbf{Second step.}
Let first compute the time-derivative of \eqref{z_2}. For that, we use the result of lemma \ref{lemma} (see equation \eqref{eq:lemme}). Then we obtain
\begin{equation} \label{dotz2}
    \dot{z}_2= u + K_1|x_1|^{|x_1|^\beta+\beta-1} (1+\beta\, ln(|x_1|)) \dot{x}_1.
\end{equation}}
\noindent \MG{Now considering the Lyapunov function 
\begin{equation} \label{v2}
 V_2=V_1+\frac{z_2^2}{2}.   
\end{equation} 
The time-derivative of $\eqref{v2}$ is:
\begin{eqnarray} \label{dotV2}
    \dot{V}_2&=& - K_1|x_1|^{1+|x_1|^\beta} + x_1 z_2  \label{eq:dot_V2} \\
    &+& z_2 \left(u + K_1|x_1|^{|x_1|^\beta+\beta-1} \right.\nonumber  \\ 
    & &\left.(1+\beta\, ln(|x_1|)) (-K_1\lceil x_1 \rfloor^{|x_1|^\beta}+z_2)\right). \nonumber
\end{eqnarray}}
\MG{By setting $u$ as proposed in \eqref{pt2} (including \eqref{ptx1}), \eqref{dotV2} becomes
\begin{equation}
    \dot{V}_2=- K_1|x_1|^{1+|x_1|^\beta} -K_2| z_2|^{1+|z_2|^\gamma} \label{eq:dot_v2_final}
\end{equation}
From \eqref{eq:dot_v2_final} the guaranteed-finite time convergence follows and this end the proof. \hfill $\Box$\\}

From Theorem \ref{Theorem}, we can set our first corollary:\\

\begin{Corollary} 
    If $d$ is bounded and its bound is know i.e. $|d| < D_{max}$,  the following control law:
    \begin{eqnarray} \label{pt2dv}
       u&=&- K_2 \sgn(z_2) - \left( \frac{}{} x_1 +  K_1|x_1|^{|x_1|^\beta+\beta-1} \right.\nonumber \\
       & &\left. (1+\beta\, ln(|x_1|)) (-K_1\lceil x_1 \rfloor^{|x_1|^\beta}+z_2) \frac{}{}\right) 
    \end{eqnarray}
    with $K_2 > D_{max}$, $\beta>1$, $K_1>0$,  and $z_2=x_2+ K_1\lceil x_1 \rfloor^{|x_1|^\beta}$ ensures a \MG{finite time} convergence of \eqref{eq:system_ref} to $x_1=x_2=0$. \vspace*{5mm}
\end{Corollary}
\textbf{Proof}: \MG{By using the same Lyapunov function $V_2$ defined in \eqref{v2}}, we obtain:
\begin{equation}
    \dot{V}_2= - K_1|x_1|^{\MG{|x_1|^\beta}} - z_2 (- K_2 \sgn(z_2) + d)
\end{equation}
and as $K> D_{max}$, we have a convergence in finite time but note in guaranteed time. \hfill $\Box$\\

\begin{Corollary} 
    If $d$ is \MG{constant, bounded and its bound is unknown}, the following control law:
    \begin{eqnarray} \label{pt2dc}
       u&=&w-K_2\lceil z_2 \rfloor^{|z_2|^\gamma}- \left(\frac{}{} x_1 + K_1|x_1|^{|x_1|^\beta+\beta-1} \right. \nonumber\\
       &&\left.  (1+\beta\, ln(|x_1|)) (-K_1\lceil x_1 \rfloor^{|x_1|^\beta}+z_2) \frac{}{}\right) 
       \\
       \dot{w}&=& - z_2 \nonumber
    \end{eqnarray}
    with $\beta>1$, $\gamma >0$ and $z_2=x_2+ K_1\lceil x_1 \rfloor^{|x_1|^\beta}$ ensures a \MG{guaranteed}-finite time convergence of \eqref{eq:system_ref} to $x_1=x_2=0$. \vspace*{5mm}
\end{Corollary}
\textbf{Proof}: \MG{By considering a new Lyapunov function as follows \begin{equation}
    V_3=V_2+\frac{(d-w)^2}{2},
\end{equation}
where $V_2$ is defined in \eqref{v2},  we get:}
\begin{equation}
    \dot{V}_3 = - K_1|x_1|^{1+|x_1|^\beta} -K_2|z_2|^{1+|z_2|^\gamma}. \label{dot_V_2_LaSalle}
\end{equation}
From LaSalle theorem, we deduce that the system \eqref{eq:system_ref}  controlled by the input \eqref{pt2dc} converges in \MG{guaranteed}-finite time to the invariant set $IS$:
\begin{equation}
IS=\{ x_1=0, \,  x_2=0\, \quad \mathbf{and}\quad  w\in \mathbb{R}\} \label{IS_corollary_1} 
\end{equation}
This ends the proof. \hfill $\Box$\\

We end this section with a remark that highlights the usefulness of the power tower function in the case of terminal sliding mode.
\begin{Remark}
 Let us consider the system \eqref{eq:system_ref}, taking the following terminal sliding surface \cite{Man94} based on a power tower function
 \begin{equation}
     s=x_2+k \lceil x_1 \rfloor^{|x_1|^\beta}   \label{TSM}
 \end{equation}
 instead of the one proposed in \cite{Man94}
 \begin{equation*}
     s=x_2+k \lceil x_1 \rfloor^{q/p}, 
 \end{equation*}
 with $p>q>0$. The time derivative of $s$ defined in \eqref{TSM} is
 \begin{equation}
     \dot{s}=u+d+k |x_1|^{|x_1|^\beta+\beta-1} (1+\beta ln(|x_1|) x_2. \label{d_TSM}
 \end{equation}
 From \eqref{d_TSM}, we can deduce following power tower terminal control law
  \begin{equation}
     u=-\lceil s \rfloor^{|s|^\gamma}-k |x_1|^{|x_1|^\beta+\beta-1} (1+\beta ln(|x_1|) x_2, \label{u_TSM}
 \end{equation}
 which has no singularity at $x_1=0$ for $\beta>1$, because $lim_{|x_1|\rightarrow 0}|x_1|^{\beta-1}  ln(|x_1|)=0$. The choice of $\beta>1$ bypasses the singularity in a different way that the one proposed in \cite{Fen02,Fen13} (i.e., define an equivalent sliding mode surface $s=x_1+k \lceil x_2 \rfloor^{p/q}$, $1<p/q<2$).
\end{Remark}

\section{Simulations} \label{simu}
\MG{To test the validity of the proposed power tower control, 4 simulations are presented. They are conducted using Matlab software, with a solver based on explicit Euler type where the sampling time is fixed to $50\mu s$. System \eqref{eq:system_ref} is considered with the following initial conditions: $x_{1}(0)=1$,
$x_{2}(0)=-1.5$. \MG{The control ``gain" is very high when the initial conditions are  far from zero, which requires very small sample step. As our sample step is limited to $10^{-6}\, sec$ and our solver is an explicit Euler scheme, we have taken initial conditions not too far from zero.}  When the perturbation $d$ is different from zero, controllers \eqref{pt2dc} and \eqref{pt2dv} are used, where $w(t=0)=0$  in \eqref{pt2dc} and $K=10$ in \eqref{pt2dv}. }
\subsection{Results with $d=0$ and with sign function}
\begin{figure}[ht!]    \centering
\hspace*{-1cm}
\includegraphics[width=1\columnwidth]{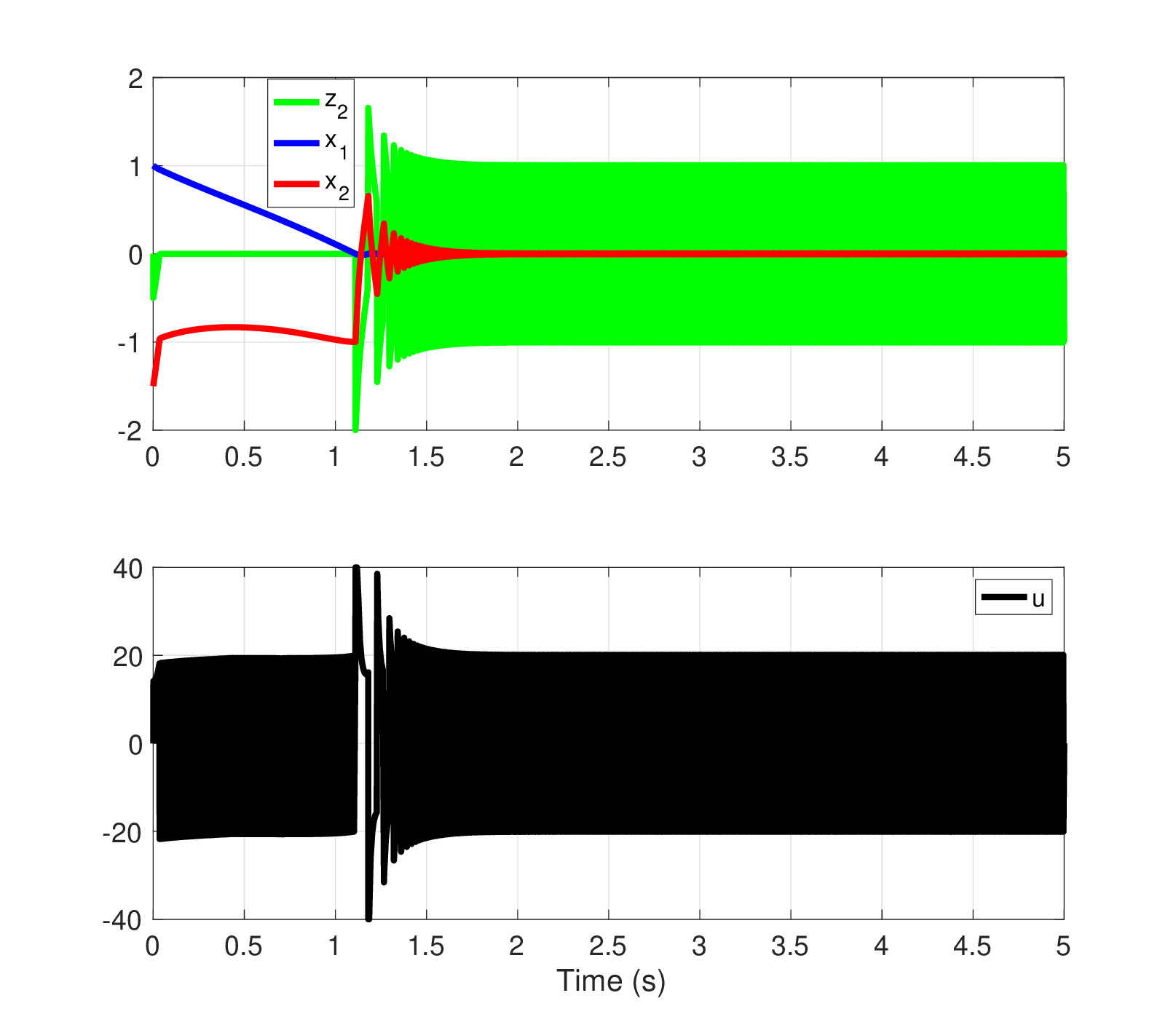}
\caption{$z_2$, $x_1$, $x_2$ (top) and $u$ (bottom)}
    \label{PT-sd-signe}
\end{figure}

In this part, the controller \eqref{pt2} is applied to \eqref{eq:system_ref}. The control gains are selected as follows $\beta=2$, 
$\gamma=1.5$, $K_1=1$ and $K_2=20$. The obtained results are depicted in Fig. \ref{PT-sd-signe}. We can notice the very good performance of the proposed controller. The state $x_1$ converges to zero in a \MG{guaranteed}-finite time. The same conclusion is stated for $x_2$, which converges to zero in \MG{fixed}-finite time. As $z_2$ is function of $x_1$ and $x_2$, its convergence to zero is also achieved. For the behavior of the control \eqref{pt2}, we can observe a chattering phenomena in $z_2$ when $x_1$ converges to zero at $t=1.2s$. This behavior is natural and can be explained by \eqref{eq:limPT_a=0}. Moreover,  at this time, when $x_1=0$, a peak on $u$ is observed, it comes from the fact that in the control law, to avoid any problems with $ln(x_1)$ in the vicinity of $x_1=0$, we cancel the product $|x_1|^{\beta-1} ln(|x_1|)$ around $x_1=0$.  To overcome the chattering behavior, we propose to replace the discontinuous function of this control by a continuous one when $z_2$ and $x_1$ reaches zero. This introduces the next subsection.    

\subsection{Results with $d=0$ and with tanh function}
\begin{figure}[ht!]    \centering
\hspace*{-1cm}
\includegraphics[width=1\columnwidth]{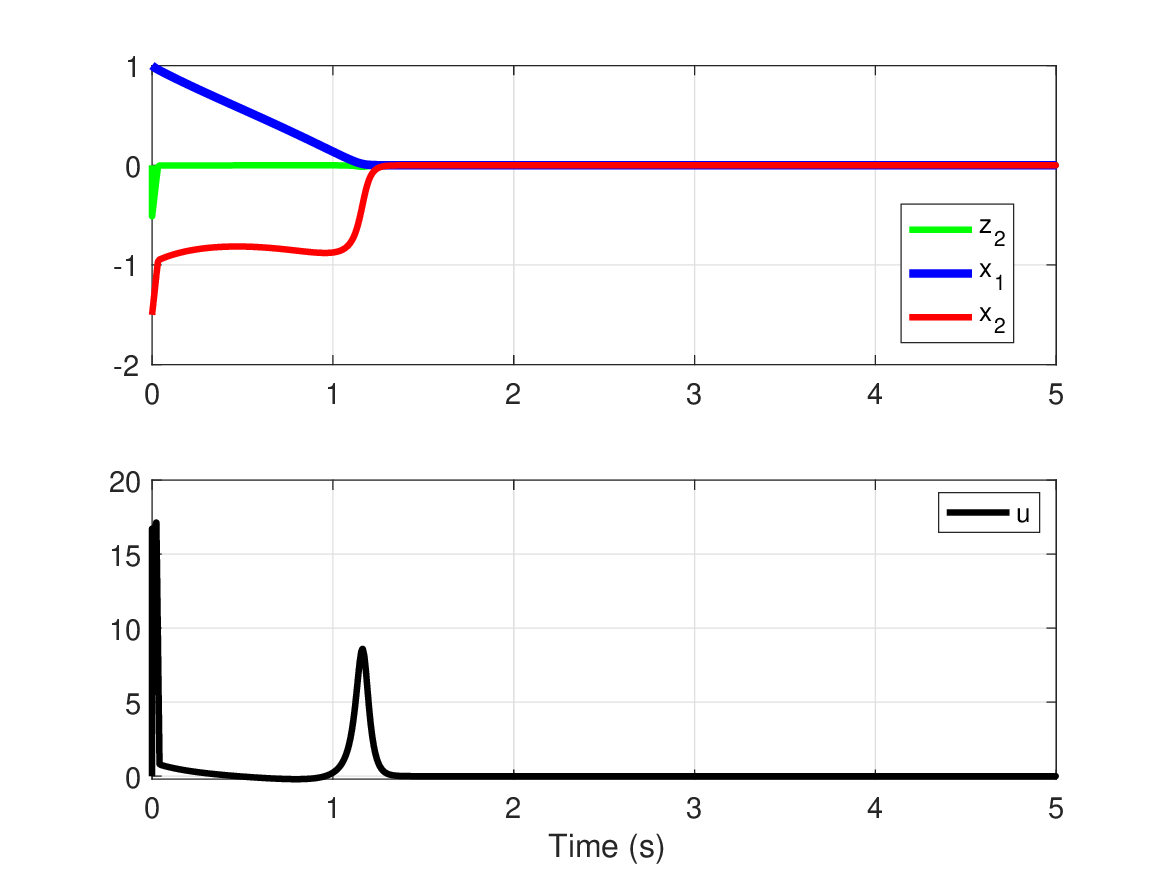}
\caption{$z_2$, $x_1$, $x_2$ (top) and $u$ (bottom)}    \label{PT-sd-tanh}
\end{figure}

\begin{figure}[ht!]    \centering
\hspace*{-1cm}
\includegraphics[width=1\columnwidth]{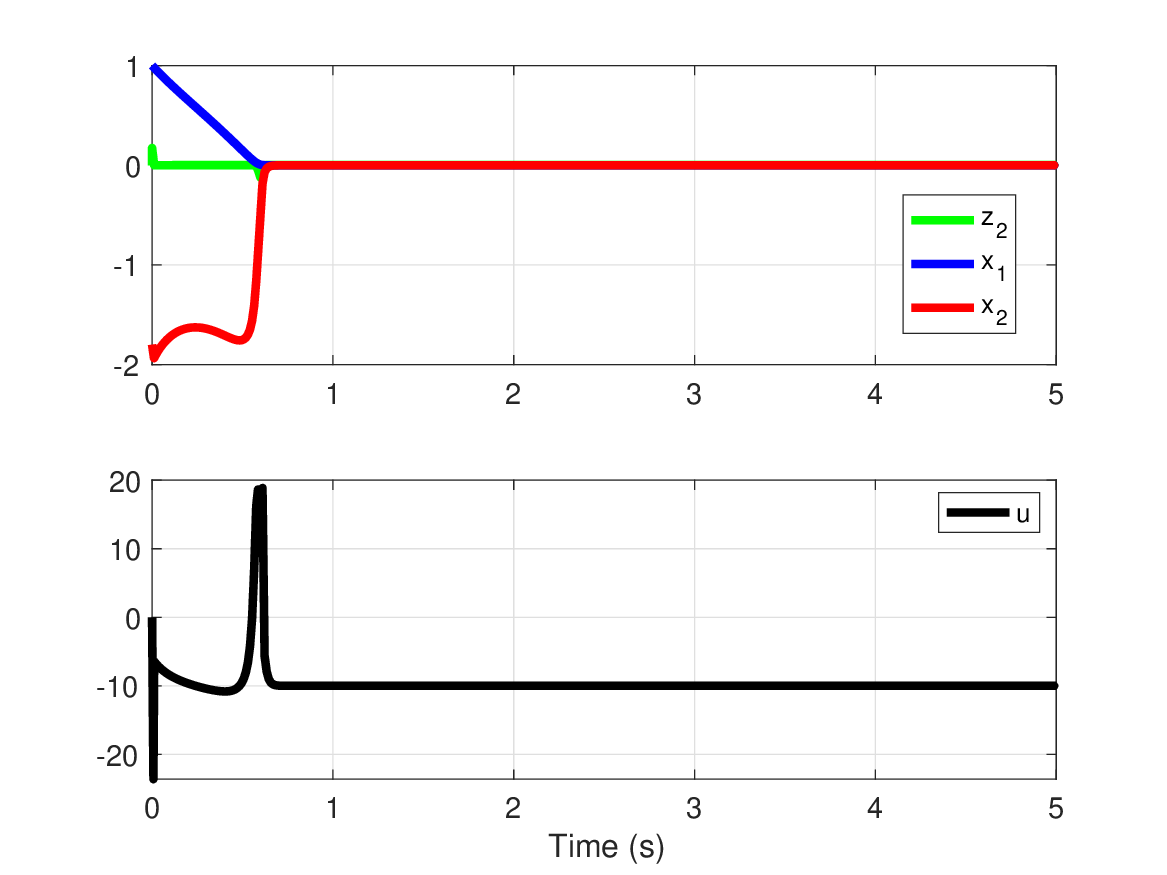}
\caption{$z_2$, $x_1$, $x_2$ (top) and $u$ (bottom)}
    \label{PT-sd-tanh-dc}
\end{figure}

\begin{figure}[ht!]    \centering
\hspace*{-1cm}
\includegraphics[width=1\columnwidth]{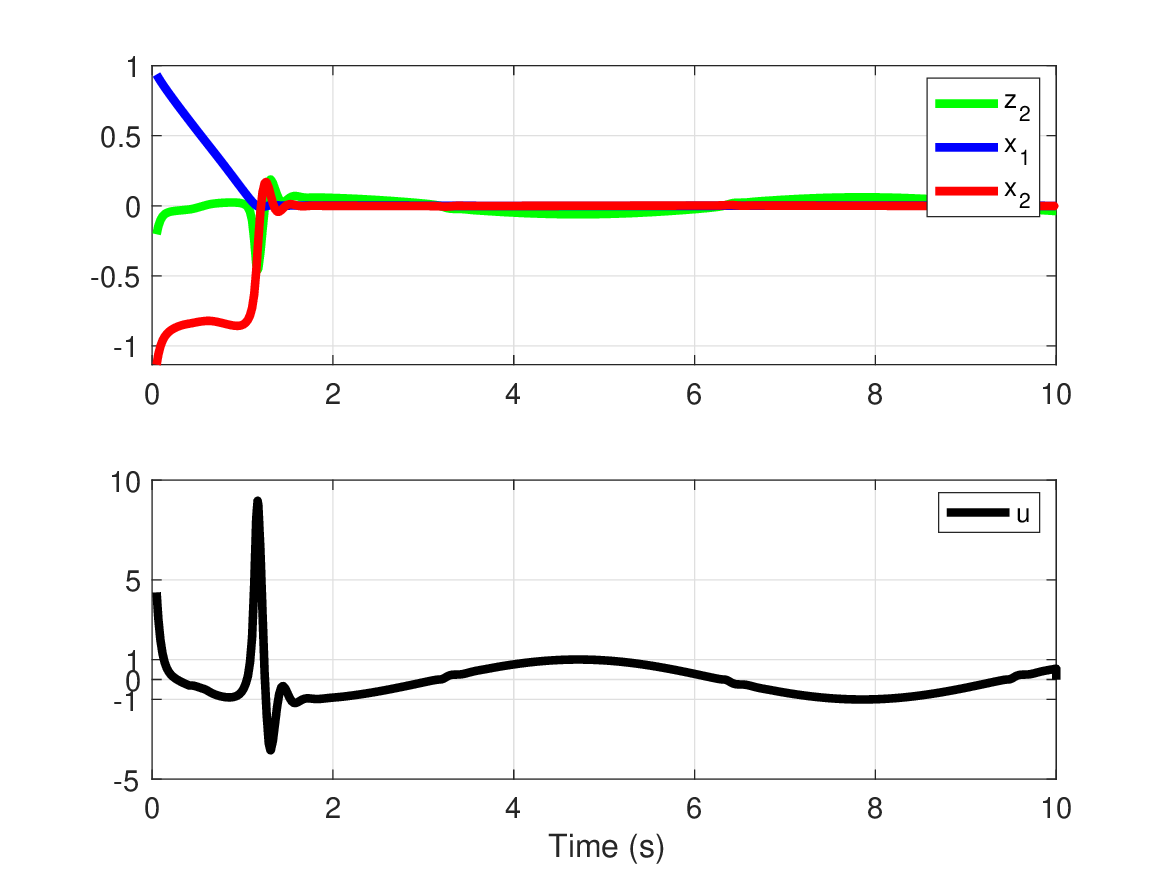}
\caption{$z_2$, $x_1$, $x_2$ (top) and $u$ (bottom)}
    \label{PT-sd-tanh-dv}
\end{figure}

\MG{In this part, the sign function used in \eqref{pt2} is replaced by the tanh function, where the gain of this function is fixed to $50$ to be more close to the behaviour of the sign function. The obtained results \JPB{a}re depicted in Fig. \ref{PT-sd-tanh}. As excepted the same results about the \MG{guaranteed}-finite time convergences of the states $x_1$, $x_2$ and $z_2$ to zero are  obtained. However, we can notice that the chattering disappeared in the control $u$, and the pic is reduced thanks to the tanh function. 
}
%
%
\MG{\subsection{Results with $d\neq0$ constant and with tanh function}}
In this part, we test the performances of the controller defined in \eqref{pt2dc} when it is applied to system \eqref{eq:system_ref} in presence of a constant bounded perturbation $d$.The control gains as selected as follows: $\beta=2$, $\gamma=1.5$, $K_1=1$ and $K_2=20$. We fixed $d=10$ (\JPB{any} another value can be chosen). The initial condition of the integrator in \eqref{pt2dc} is chosen equal to zero ($w_0=0$). The simulation results are shown in \eqref{PT-sd-tanh-dc}. In Fig. \ref{PT-sd-tanh-dc} we can show that the control \eqref{pt2dc} performs well in the sense that the perturbation $d$ is exactly canceled ($u=-d$ in steady state) thanks to the integral term $w$ in \eqref{pt2dc}. The latter is replaced by the tanh function with high gain (50) to avoid chattering phenomenon. 
Even if the $sign$ function is approximated by the $tanh$ function for avoiding the chattering, the convergence of the states $x_1$ and $x_2$ seems to be in guaranteed-finite time. 
%
%
\MG{\subsection{Results with $d$ variable and with tanh function}}
\MG{In this part, we test the performances of the controller defined in \eqref{pt2dv} when it is applied to system \eqref{eq:system_ref} in presence of a variable bounded perturbation $d$ where its upper bound is known. For that we took $d$ as a sinus function: $d=sin(t)$ and $\beta=2$ and $\gamma=0$. The control gains are fixed $\beta=2$, $\gamma=0$, $K_1=1$ and $K_2=20$. Then, $K_2>D_{max}$, where $D_{max}=1$ is the upper bound of the sinus function.   The obtained results, plotted in Fig. \eqref{PT-sd-tanh-dv},
  show very good performances of the control \eqref{pt2dv}, the perturbation $d$ is exactly canceled ($u=-d$) thanks to the sign function in \eqref{pt2dc}, which is replaced by tanh function (gain$=50$) to avoid chattering phenomenon.  In this case the convergence of the states $x_1$ and $x_2$ are in finite time and not in guaranteed time \MG{for the $sign$ function and only asymptotic for $tanh$ function. It can be noticed also that $z_2$ has not converged exactly to zero. This is due to the fact the gain of the $tanh$ function is not chosen so big to avoid control pics.    }}

\section{Conclusion} \label{conc}
In this paper, we propose a new control law on a power tower function truncated to order 2. This function makes it possible to use the backstepping technique in order to propose convergence in guaranteed and finite time. 
In a future work, the convergence finite times  including the fixed one will be derived and the non-matching perturbations problem will be studied. The fixed-time will be computed in function of $K_1$ and $K_2$ which multiply respectively the term $\lceil x_1 \rfloor^{|x_1|^\beta}$ and the term $\lceil z_2 \rfloor^{|z_2|^\gamma}$ in the control design. Moreover, the approach will be extended to higher dimensional systems.  In a second step, the behavior of such law in an observer-based control scheme or with respect to noisy measurements or actuator saturation or again under sampling must be investigated.

\bibliographystyle{plain}
\bibliography{Biblio}
\end{document}